\documentclass[10pt]{article}
\usepackage{amsmath, amssymb, amsthm}

\newtheorem{theorem}{Theorem}

\theoremstyle{definition}

\title{A combinatorial proof of the quadratic Gauss sum}
\author{Adan Medrano Martin del Campo}
\date{September 29, 2026}

\begin{document}

\maketitle

\begin{abstract}
We provide an elementary approach to computing the quadratic Gauss sum, emphasizing the method employed to compute the sum corresponding to multiples of $4$. The proof method employed does not rely on the law of quadratic reciprocity, but only on exploiting symmetries among the sum's components, together with a recursive approach.
\end{abstract}

\vspace{10px}


\vspace{10px}

Let $n$ be a positive integer and $\zeta_{n} = e^{\frac{2\pi i}{n}}$ be a primitive $n$-th root of unity. The \emph{quadratic Gauss sum} is defined as
\[
    \mathcal{G}(n) = \sum_{k = 0}^{n - 1} \zeta_{n}^{k^{2}}.
\]
The calculation of the value of $\mathcal{G}(n)$ for every positive integer $n$ is a classical result of Gauss, summarized in the following theorem.

\begin{theorem}{(Gauss)}
    Let $n$ be a positive integer. Then
    \[
        \frac{1}{\sqrt{n}}\mathcal{G}(n) =
        \begin{cases}
            1     & n\equiv 1\pmod{4}\\
            0     & n\equiv 2\pmod{4}\\
            i     & n\equiv 3\pmod{4}\\
            1 + i & n\equiv 0\pmod{4}
        \end{cases}
    \]
\end{theorem}

Many proofs of this result exist. We present a classical linear-algebraic proof studied by many authors (\cite{Schur}, \cite{Waterhouse}, \cite{Murty}, \cite{MurtyPathak}) for $n$ not divisible by $4$. For $n$ divisible by $4$, we present a proof that relies only on symmetry, combinatorics, and recursion. 

\section{Fourier-type matrices}

    For $\zeta\in \mathbb{C}$, consider the $n\times n$ matrix
    \[
        \mathcal{F}_{n}(\zeta) = \frac{1}{\sqrt{n}}
        \begin{pmatrix}
            1      & 1             & 1                & \cdots & 1 \\
            1      & \zeta         & \zeta^{2}        & \cdots & \zeta^{n - 1} \\
            1      & \zeta^{2}     & \zeta^{4}        & \cdots & \zeta^{2(n - 1)} \\
            \vdots & \vdots        & \vdots           & \ddots & \vdots \\
            1      & \zeta^{n - 1} & \zeta^{2(n - 1)} & \cdots & \zeta^{(n - 1)^{2}}
        \end{pmatrix}
    \]
    where the $0$-indexed entry of this matrix in its $a$-th row and $b$-th column is $\frac{1}{\sqrt{n}}\zeta^{ab}$. Letting $\zeta = \zeta_{n}$ gives the \emph{discrete Fourier transform} matrix, studied for example in \cite{Fourier}, and it is clear that
    \[
        \frac{1}{\sqrt{n}}\mathcal{G}(n) = \text{Tr}(\mathcal{F}_{n}(\zeta_{n})).
    \]
    We may calculate this trace using only linear-algebraic methods for values of $n$ not divisible by $4$. For this, we study matrices $\mathcal{F}_{n}(\zeta_{n}^{m})$ with $\gcd(n, m) = 1$.

    \subsection{Eigenvalues of $\mathcal{F}_{n}(\zeta_{n}^{m})$}

        The dot product between row $a$ and column $b$ of $\mathcal{F}_{n}(\zeta_{n}^{m})$ is
        \[
            \frac{1}{n}\left(\sum_{k = 0}^{n - 1} \zeta_{n}^{(a + b)km}\right) =
            \begin{cases}
                1 & n\mid  a + b \\
                0 & n\nmid a + b 
            \end{cases}
        \]
        hence, the square of our matrix is the permutation matrix,
        \[
            \mathcal{F}_{n}(\zeta_{n}^{m})^{2} =
            \begin{pmatrix}
                1      & 0      & 0      & \cdots & 0      & 0 \\
                0      & 0      & 0      & \cdots & 0      & 1 \\
                0      & 0      & 0      & \cdots & 1      & 0 \\
                \vdots & \vdots & \vdots & \ddots & \vdots & \vdots \\
                0      & 0      & 1      & \cdots & 0      & 0 \\
                0      & 1      & 0      & \cdots & 0      & 0 \\
            \end{pmatrix}
        \]
        which consists of $2$-cycles and fixed points, and thus $\mathcal{F}_{n}(\zeta_{n}^{m})^{4}$ is the identity. Hence, the eigenvalues of $\mathcal{F}_{n}(\zeta_{n}^{m})$ are in the set $\{\pm 1, \pm i\}$ and the eigenvalues of $\mathcal{F}_{n}(\zeta_{n}^{m})^{2}$ must be in the set $\{\pm 1\}$. In particular, it is clear that
        \[
            \text{Tr}(\mathcal{F}_{n}(\zeta_{n}^{m})^{2}) = 
            \begin{cases}
                1 & n\equiv 1\pmod{2} \\
                2 & n\equiv 0\pmod{2}
            \end{cases}
        \]
    \newpage
    \subsection{Trace norm of $\mathcal{F}_{n}(\zeta_{n}^{m})$}

        We can compute the norm of $\mathcal{F}_{n}(\zeta_{n}^{m})$ by observing
        \begin{align*}
            \left|\sum_{k = 0}^{n - 1}\zeta_{n}^{mk^{2}}\right|^{2} &=
            \left(\sum_{a = 0}^{n - 1}\zeta_{n}^{ma^{2}}\right)
            \left(\sum_{b = 0}^{n - 1}\zeta_{n}^{-mb^{2}}\right) \\ &= 
            \sum_{a = 0}^{n - 1}\sum_{b = 0}^{n - 1}\zeta_{n}^{m(a + b)(a - b)} \\ &=
            \sum_{d = 0}^{n - 1}\sum_{b=0}^{n - 1}\zeta_{n}^{md(d + 2b)} \\ &=
            \sum_{d = 0}^{n - 1}\zeta_{n}^{md^{2}}\left(\sum_{b=0}^{n - 1}\zeta_{n}^{2mdb}\right) \\ &=
            \sum_{n\mid 2d} \zeta_{n}^{md^{2}}n =
            \begin{cases}
                n                       & n\equiv 1\pmod{2} \\
                n + (-1)^{\frac{n}{2}}n & n\equiv 0 \pmod{2}
            \end{cases}
        \end{align*}
        Hence,
        \[
            \left|\text{Tr}(\mathcal{F}_{n}(\zeta_{n}^{m})\right|^{2} =
            \begin{cases}
                0 & n \equiv 2\pmod{4} \\
                1 & n \equiv 1\pmod{2} \\
                2 & n \equiv 0\pmod{4}
            \end{cases}
        \]

    \subsection{Eigenspaces of $\mathcal{F}_{n}(\zeta_{n}^{m})$}
    
        Since
        \[
            \mathcal{F}_{n}(\zeta_{n}^{m})^{-1} =
            \mathcal{F}_{n}(\zeta_{n}^{-m}) =
            \mathcal{F}_{n}(\zeta_{n}^{m})^{\dagger}
        \]
        it follows that $\mathcal{F}_{n}(\zeta_{n}^{m})$ is unitary, hence diagonalizable. We can then study the eigenvalue multiplicity of $\mathcal{F}_{n}(\zeta_{n}^{m})$ by letting it act on a $\mathbb{C}$-vector space $V_{n}$, and studying its eigenspace decomposition
        \[
            V_{n} \cong V_{n}^{1}\oplus V_{n}^{-1}\oplus V_{n}^{i} \oplus V_{n}^{-i}
        \]
        where $\mathcal{F}_{n}(\zeta_{n}^{m})$ acts on $V_{n}^{i^{k}}$ by multiplication by $i^{k}$. The matrix $\mathcal{F}_{n}(\zeta_{n}^{m})^{2}$ acts on $V_{n}^{1}\oplus V_{n}^{-1}$ as the identity and on $V_{n}^{i}\oplus V_{n}^{-i}$ as multiplication by $-1$. Letting
        \[
            \dim(V_{n}^{i^{k}}) = d_{i^{k}}
        \]
        we have $n = d_{1} + d_{-1} + d_{i} + d_{-i}$ and
        \begin{align*}
            \text{Tr}(\mathcal{F}_{n}(\zeta_{n}^{m})) &= (d_{1} - d_{-1}) + \left(d_{i}
            - d_{-i}\right)i \\
            \text{Tr}(\mathcal{F}_{n}(\zeta_{n}^{m})^{2}) &= (d_{1} + d_{-1}) - (d_{i} + d_{-i})
        \end{align*}
        Moreover, $|\text{Tr}(\mathcal{F}_{n}(\zeta_{n}^{m}))|^{2} = (d_{1} - d_{-1})^{2} + (d_{i} - d_{-i})^{2} \leq 2$, which implies that each of $|d_{1} - d_{-1}|, |d_{i} - d_{-i}|$ is at most $1$, and equal to $0$ if and only if it is even.
    
    \newpage
    \subsection{Determinant of $\mathcal{F}_{n}(\zeta_{n})$} 

        The matrix $\sqrt{n}\mathcal{F}_{n}(\zeta_{n}^{m})$ is a \emph{Vandermonde matrix} whose determinant is
        \[
            \det(\sqrt{n}\mathcal{F}_{n}(\zeta_{n}^{m})) = \prod_{0\leq a < b < n} \left(\zeta_{n}^{mb} - \zeta_{n}^{ma}\right).
        \]
        We begin by computing this determinant for $m = 1$ by manipulating each factor,
        \begin{align*}
            \zeta_{n}^{b} - \zeta_{n}^{a} &= 
            e^{\frac{2\pi i}{n}\left(\frac{a + b}{2}\right)}\left(e^{\frac{2\pi i}{n}\left(\frac{b - a}{2}\right)} - e^{-\frac{2\pi i}{n}\left(\frac{b - a}{2}\right)}\right) \\ &=
            e^{\frac{\pi i}{2}\cdot \left(\frac{2(a + b)}{n} + 1\right)}\cdot 2\sin\left(\left(\frac{b - a}{n}\right)\pi\right).
        \end{align*}
        Since $0 < \left(\frac{b - a}{n}\right)\pi < \pi$, we have $\sin\left(\left(\frac{b - a}{n}\right)\pi\right) > 0$ and hence
        \[
            \text{Arg}(\zeta_{n}^{b} - \zeta_{n}^{a}) = \frac{\pi}{2}\cdot \left(\frac{2(a + b)}{n} + 1\right) \pmod{2\pi}.
        \]
        The sum of these factors' arguments yields $\text{Arg}(\det(\sqrt{n}\mathcal{F}_{n}(\zeta_{n})))$, which is
        \begin{align*}
            \text{Arg}(\det(\sqrt{n}\mathcal{F}_{n}(\zeta_{n}))) &=  
            \sum_{0\leq a < b < n} \frac{\pi}{2}\cdot \left(\frac{2(a + b)}{n} + 1\right) \pmod{2\pi}\\ &=
            \frac{\pi}{2}\cdot \left(\binom{n}{2} + (n - 1)\left(\sum_{k = 0}^{n - 1}\frac{2k}{n}\right)\right) \pmod{2\pi}\\ &=
            \frac{\pi}{2}\cdot\left(\binom{n}{2} + (n - 1)^{2}\right) \pmod{2\pi}.
        \end{align*}
        Furthermore, $\det(\mathcal{F}_{n}(\zeta_{n}))$ has unit norm as $\mathcal{F}_{n}(\zeta_{n})$ is unitary, thus
        \[
            \det(\mathcal{F}_{n}(\zeta_{n})) = 
            \frac{1}{\sqrt{n^{n}}}\det(\sqrt{n}\mathcal{F}_{n}(\zeta_{n})) = i^{\binom{n}{2} + (n - 1)^{2}}.
        \]
        
    \subsection{Periodicity of $\det\left(\mathcal{F}_{n}(\zeta_{n})\right)$}
        Viewing $\det(\mathcal{F}_{n}(\zeta_{n}))$ as the product of the eigenvalues of $\mathcal{F}_{n}(\zeta_{n})$, we have
        \[
            \det(\mathcal{F}_{n}(\zeta_{n})) = \prod_{k = 0}^{3}(i^{k})^{d_{i^{k}}} = i^{2d_{-1} + d_{i} - d_{-i}}
        \]
        and from the computation of $\det\left(\mathcal{F}_{n}(\zeta_{n})\right)$ it follows that
        \[
            \dbinom{n}{2} + (n - 1)^{2} \equiv 2d_{-1} + d_{i} - d_{-i} \pmod{4}.
        \]
        In particular, $\binom{n}{2} + (n - 1)^{2} \pmod{4}$ is $8$-periodic since
        \begin{align*}
            \left(\dbinom{n + 8}{2} + (n + 7)^{2}\right) - \left(\dbinom{n}{2} + (n - 1)^{2}\right) &=
            4(6n + 19)
        \end{align*}
        Therefore, the sequence $(\det(\mathcal{F}_{n}(\zeta_{n})))_{n\geq 1}$ is $8$-periodic.

    \newpage
\section{Computing $\mathcal{G}(n)$ for $n$ not divisible by $4$}

        We use the dimensions of the eigenspaces of $\mathcal{F}_{n}(\zeta_{n})$ to compute its trace when $n$ is not divisible by $4$, keeping the information for $n$ viewed by $\pmod{8}$. We begin by computing $d_{1} + d_{-1}$ and $d_{i} + d_{-i}$ since
        \begin{align*}
            d_{1} + d_{-1} &= \frac{1}{2}\left(n + \text{Tr}(\mathcal{F}_{n}(\zeta_{n})^{2})\right) \\
            d_{i} + d_{-i} &= \frac{1}{2}\left(n - \text{Tr}(\mathcal{F}_{n}(\zeta_{n})^{2})\right)
        \end{align*}
        Furthermore, if $d_{1} + d_{-1}$ is even, $d_{1} = d_{-1}$ and if $d_{i} + d_{-i}$ is even, $d_{i} = d_{-i}$. Together, this can be summarized in the following table
        \begin{center}
            \begin{tabular}{| l | l | l | l | l |}
                \hline \rule{0pt}{2.5ex}
                $n$ & $d_{1}$ & $d_{-1}$ & $d_{i}$ & $d_{-i}$ \\ [0.5ex]
                \hline
                $8N + 1$ &&& $2N$ & $2N$ \\
                $8N + 2$ & $2N + 1$ & $2N + 1$ & $2N$ & $2N$ \\
                $8N + 3$ & $2N + 1$ & $2N + 1$ && \\
                \hline
                $8N + 5$ &&& $2N + 1$ & $2N + 1$ \\
                $8N + 6$ & $2N + 2$ & $2N + 2$ & $2N + 1$ & $2N + 1$ \\
                $8N + 7$ & $2N + 2$ & $2N + 2$ && \\
                \hline
            \end{tabular}
        \end{center}
        Finally, using the condition
        \[
            \dbinom{n}{2} + (n - 1)^{2} \equiv 2d_{-1} + d_{i} - d_{-i} \pmod{4}
        \]
        allows us to complete the missing dimensions for odd $n$,
        \begin{center}
            \begin{tabular}{| l | l | l | l | l |}
                \hline \rule{0pt}{2.5ex}
                $n$ & $d_{1}$ & $d_{-1}$ & $d_{i}$ & $d_{-i}$ \\ [0.5ex]
                \hline
                $8N + 1$ & $2N + 1$ & $2N$     & $2N$ & $2N$ \\
                $8N + 2$ & $2N + 1$ & $2N + 1$ & $2N$ & $2N$ \\
                $8N + 3$ & $2N + 1$ & $2N + 1$ & $2N + 1$ & $2N$ \\
                \hline
                $8N + 5$ & $2N + 2$ & $2N + 1$ & $2N + 1$ & $2N + 1$ \\
                $8N + 6$ & $2N + 2$ & $2N + 2$ & $2N + 1$ & $2N + 1$ \\
                $8N + 7$ & $2N + 2$ & $2N + 2$ & $2N + 2$ & $2N + 1$ \\
                \hline
            \end{tabular}
        \end{center}
        This proves that
        \begin{align*}
            \text{Tr}(\mathcal{F}_{4N + 1}(\zeta_{4N + 1})) &= 1 &\implies&& \mathcal{G}(4N + 1) &= \sqrt{4N + 1}\\
            \text{Tr}(\mathcal{F}_{4N + 2}(\zeta_{4N + 2})) &= 0 &\implies&& \mathcal{G}(4N + 2) &= 0\\
            \text{Tr}(\mathcal{F}_{4N + 3}(\zeta_{4N + 3})) &= i &\implies&& \mathcal{G}(4N + 3) &= \sqrt{4N + 3}\cdot i
        \end{align*}
        It only remains to compute the value of $\mathcal{G}(4N)$.

        \newpage
\section{The special case of $n$ odd and $m = 2$}
    \subsection{Determinant of $\mathcal{F}_{n}(\zeta_{n}^{2})$}

        We compute
        \[
            \det(\sqrt{n}\mathcal{F}_{n}(\zeta_{n}^{2})) = \prod_{0\leq a < b < n} \left(\zeta_{n}^{2b} - \zeta_{n}^{2a}\right)
        \]
        by counting the factor sign changes with respect to $\det(\sqrt{n}\mathcal{F}_{n}(\zeta_{n}))$. For this, we consider the following combinatorial interpretation: We list the powers
        \[
            1, \zeta_{n}^{m}, \zeta_{n}^{2m}, \ldots, \zeta_{n}^{(n - 1)m}.
        \]
        Then each factor of $\det(\sqrt{n}\mathcal{F}_{n}(\zeta_{n}^{m}))$ is given by taking every pair of distinct elements on the above list, and taking the difference between the rightmost element and the leftmost element. For $m = 2$ this list is
        \[
            1, \zeta_{n}^{2}, \ldots, \zeta_{n}^{n - 1}, \zeta_{n}, \ldots, \zeta_{n}^{n - 2}
        \]
        so the element $1$ is fixed, and the odd powers and even powers of $\zeta_{n}$ preserve their ordering. Hence, the only pairs that have a sign change consist of an odd and an even power of $\zeta_{n}$. It suffices to count, for each even power $\zeta_{n}^{2k}$ with $0 < k \leq \frac{n - 1}{2}$, how many odd powers $\zeta_{n}^{2l - 1}$ with $0 < l \leq \frac{n - 1}{2}$ and $2l - 1 < 2k$ lie to its right in the list. This is equivalent to $0 < l \leq k$ and hence there are exactly $k$ such pairs for each $k$. Hence, the total number of sign flips is
        \[
            \sum_{k = 1}^{\frac{n - 1}{2}} k = \dbinom{\frac{n + 1}{2}}{2} = \frac{n^{2} - 1}{8}
        \]
        and therefore, for odd $n$ we have
        \[
            \det(\mathcal{F}_{n}(\zeta_{n}^{2})) = (-1)^{\frac{n^{2} - 1}{8}}\det(\mathcal{F}_{n}(\zeta_{n})).
        \]

    \subsection{Trace of $\mathcal{F}_{n}(\zeta_{n}^{2})$} 

        We can compute the dimensions of the eigenspaces of $\mathcal{F}_{n}(\zeta_{n}^{2})$ for odd $n$ analogously to those of $\mathcal{F}_{n}(\zeta_{n})$. Using the dimensions computed for $\mathcal{F}_{n}(\zeta_{n})$ together with $\det(\mathcal{F}_{n}(\zeta_{n}^{2}))$ expressed in terms of $\det(\mathcal{F}_{n}(\zeta_{n}))$, we obtain the following
        \begin{center}
            \begin{tabular}{| l | l | l | l | l |}
                \hline \rule{0pt}{2.5ex}
                $n$ & $d_{1}$ & $d_{-1}$ & $d_{i}$ & $d_{-i}$ \\ [0.5ex]
                \hline
                $8N + 1$ & $2N + 1$ & $2N$     & $2N$     & $2N$ \\
                $8N + 3$ & $2N + 1$ & $2N + 1$ & $2N$     & $2N + 1$ \\
                $8N + 5$ & $2N + 1$ & $2N + 2$ & $2N + 1$ & $2N + 1$ \\
                $8N + 7$ & $2N + 2$ & $2N + 2$ & $2N + 2$ & $2N + 1$ \\
                \hline
            \end{tabular}
        \end{center}
        Hence, for odd $n$ we have
        \[
            \text{Tr}(\mathcal{F}_{n}(\zeta_{n}^{2})) = (-1)^{\frac{n^{2} - 1}{8}}\text{Tr}(\mathcal{F}_{n}(\zeta_{n})).
        \]

\newpage
\section{Recursion for $\mathcal{G}(4N)$}

    We steer away from linear algebra now and show the following recurrence
    \[
        \mathcal{G}(4N) = 
        \begin{cases}
            2\mathcal{G}(N) & N\equiv 0 \pmod{4} \\
            2\mathcal{G}(N) + 2i^{N}\mathcal{G}(N) & N\equiv 1 \pmod{2} \\
            4\zeta_{8}^{2\left(\frac{N + 2}{4}\right)^{2} - 1}\mathcal{G}\left(\frac{N}{2}\right) & N\equiv 2\pmod{4}
        \end{cases}
    \]
    which allows us to prove inductively that
    \[
        \mathcal{G}(4N) = \sqrt{4N}(1 + i).
    \]
    We separate the components of $\mathcal{G}(4N)$ by the parity of their square exponent,
    \begin{align*}
        \mathcal{G}(4N) &= \sum_{k = 0}^{4N - 1}\zeta_{4N}^{k^{2}} \\ &=
        \sum_{k = 0}^{2N - 1}\zeta_{4N}^{(2k)^{2}} + \sum_{k = 0}^{2N - 1} \zeta_{4N}^{(2k + 1)^{2}} \\ &=
        \sum_{k = 0}^{2N - 1}\zeta_{N}^{k^{2}} + \sum_{k = 0}^{2N - 1} \zeta_{4N}^{(2k + 1)^{2}} \\ &=
        2\sum_{k = 0}^{N - 1}\zeta_{N}^{k^{2}} + \sum_{k = 0}^{2N - 1} \zeta_{4N}^{(2k + 1)^{2}} =
        2\mathcal{G}(N) + \sum_{k = 0}^{2N - 1} \zeta_{4N}^{(2k + 1)^{2}}.
    \end{align*}
    The key is then to manipulate the sum consisting of odd square exponents.

    \begin{itemize}
        \item If $N$ is odd, we shift exponents by $\frac{N - 1}{2}$ and obtain
            \[
                \zeta_{4N}^{(2k + N)^{2}} = 
                \zeta_{N}^{k^{2}}\cdot i^{N}
            \]
            so it follows that
            \begin{align*}
                \sum_{k = 0}^{2N - 1}\zeta_{4N}^{(2k + 1)^{2}} = 
                \sum_{k = 0}^{2N - 1}\zeta_{4N}^{(2k + N)^{2}} = 
                \sum_{k = 0}^{2N - 1}\zeta_{N}^{k^{2}}\cdot i^{N} = 
                2i^{N}\left(\sum_{k = 0}^{N - 1}\zeta_{N}^{k^{2}}\right)
            \end{align*}
            and therefore $\mathcal{G}(4N) = 2\mathcal{G}(N) + 2i^{N}\mathcal{G}(N)$.
    
        \item If $N$ is divisible by $4$, shifting exponents by $\frac{N}{2}$ we obtain
            \begin{align*}
                \zeta_{4N}^{(2k + 1 + N)^{2}} = 
                \zeta_{4N}^{(2k + 1)^{2}}\cdot (-1)^{2k + 1}\cdot \zeta_{4}^{N} = -\zeta_{4N}^{(2k + 1)^{2}}.
            \end{align*}
            Hence,
            \begin{align*}
                \sum_{k = 0}^{2N - 1}\zeta_{4N}^{(2k + 1)^{2}} = 
                \frac{1}{2}\left(\sum_{k = 0}^{2N - 1}\zeta_{4N}^{(2k + 1)^{2}} + \zeta_{4N}^{(2k + 1 + N)^{2}}\right) = 
                0
            \end{align*}
            and therefore $\mathcal{G}(4N) = 2\mathcal{G}(N)$.
    
        \item If $N\equiv 2\pmod{4}$ then $\mathcal{G}(N) = 0$. Shifting exponents by $\frac{N}{2}$ we obtain
            \begin{align*}
                \zeta_{4N}^{(2k + 1 + N)^{2}} = 
                \zeta_{4N}^{(2k + 1)^{2}}\cdot (-1)^{2k + 1}\cdot \zeta_{4}^{N} = \zeta_{4N}^{(2k + 1)^{2}}.
            \end{align*}
            Hence, the terms in the sum are $\frac{N}{2}$-periodic. In particular,
            \begin{align*}
                \sum_{k = 0}^{2N - 1}\zeta_{4N}^{(2k + 1)^{2}} = 2\sum_{k = 0}^{N - 1}\zeta_{4N}^{(2k + 1)^{2}}.
            \end{align*}
            We now make an additional exponent shift by $\frac{N - 2}{4}$ to see that
            \begin{align*}
                \zeta_{4N}^{(2k + \frac{N}{2})^{2}} = 
                \zeta_{N}^{k^{2}}\cdot (-1)^{k}\cdot \zeta_{8}^{\frac{N}{2}}
            \end{align*}
            and we have an additional symmetry as
            \[
                (-1)^{k + \frac{N}{2}}\zeta_{N}^{(k + \frac{N}{2})^{2}} = 
                \zeta_{N}^{k^{2}}\cdot (-1)^{k + N} = 
                (-1)^{k}\zeta_{N}^{k^{2}}
            \]
            so we can exclusively pick whichever of $k, k + \frac{N}{2}$ is even to obtain
            \begin{align*}
                2\sum_{k = 0}^{N - 1}\zeta_{4N}^{(2k + 1)^{2}} &=
                2\zeta_{8}^{\frac{N}{2}}\left(\sum_{k = 0}^{N - 1}\zeta_{N}^{k^{2}}\cdot(-1)^{k}\right) \\ &=
                4\zeta_{8}^{\frac{N}{2}}\left(\sum_{k = 0}^{\frac{N}{2} - 1}\zeta_{N}^{(2k)^{2}}\right) = 
                4\zeta_{8}^{\frac{N}{2}}\left(\sum_{k = 0}^{\frac{N}{2} - 1}\zeta_{\frac{N}{2}}^{2k^{2}}\right).
            \end{align*}
            Since $\text{Tr}(\mathcal{F}_{\frac{N}{2}}(\zeta_{\frac{N}{2}}^{2})) = (-1)^{\frac{N^{2} - 4}{32}}\text{Tr}(\mathcal{F}_{\frac{N}{2}}(\zeta_{\frac{N}{2}}))$ as $\frac{N}{2}$ is odd, we have
            \[
                \left(\sum_{k = 0}^{\frac{N}{2} - 1}\zeta_{\frac{N}{2}}^{2k^{2}}\right) = (-1)^{\frac{N^{2} - 4}{32}}\left(\sum_{k = 0}^{\frac{N}{2} - 1}\zeta_{\frac{N}{2}}^{k^{2}}\right).
            \]
            Putting everything together, we conclude that
            \begin{align*}
                \mathcal{G}(4N) &= 2\mathcal{G}(N) + \sum_{k = 0}^{2N - 1}\zeta_{4N}^{(2k + 1)^{2}} \\ &=
                2\sum_{k = 0}^{N - 1}\zeta_{4N}^{(2k + 1)^{2}} \\ &=
                4\zeta_{8}^{\frac{N}{2}}\left(\sum_{k = 0}^{\frac{N}{2} - 1}\zeta_{\frac{N}{2}}^{2k^{2}}\right) \\ &=
                4\zeta_{8}^{\frac{N}{2}}(-1)^{\frac{N^{2} - 4}{32}}\left(\sum_{k = 0}^{\frac{N}{2} - 1}\zeta_{\frac{N}{2}}^{k^{2}}\right) = 4\zeta_{8}^{2\left(\frac{N + 2}{4}\right)^{2} - 1}\mathcal{G}\left(\frac{N}{2}\right).
            \end{align*}
    \end{itemize}
    
\newpage
\section{Computing $\mathcal{G}(4N)$}

    We proceed inductively over $\nu_{2}(N)$. If $\nu_{2}(N) \geq 2$,
    \begin{align*}
        \mathcal{G}(4N) = 
        2\mathcal{G}(N) = 
        2\cdot \sqrt{N}(1 + i) = 
        \sqrt{4N}(1 + i).
    \end{align*}
    Hence, it suffices to prove the base cases $\nu_{2}(N) = 0, 1$.

    \begin{itemize}
        \item If $\nu_{2}(N) = 0$ we distinguish two cases
            \[
                \mathcal{G}(4N) = 
                2(1 + i^{N})\mathcal{G}(N) = 
                \begin{cases}
                    2(1 + i)\sqrt{N}        & N\equiv 1\pmod{4} \\
                    2(1 - i)\sqrt{N}\cdot i & N\equiv 3\pmod{4}
                \end{cases}
            \]
            
        \item If $\nu_{2}(N) = 1$ we distinguish two cases
            \[
                \mathcal{G}(4N) = 
                4\zeta_{8}^{2\left(\frac{N + 2}{4}\right)^{2} - 1}\mathcal{G}\left(\frac{N}{2}\right) =
                \begin{cases}
                    \sqrt{8}(1 + i)\sqrt{\frac{N}{2}}        & \frac{N}{2}\equiv 1\pmod{4} \\
                    \sqrt{8}(1 - i)\sqrt{\frac{N}{2}}\cdot i & \frac{N}{2}\equiv 3\pmod{4}
                \end{cases}
            \]
    \end{itemize}
    All four cases yield $\mathcal{G}(4N) = \sqrt{4N}(1 + i)$, thus completing the induction.

\section*{Acknowledgments} The author thanks Carlos Andres Cortez Lemos and Eli Sadovnik for useful discussions that sparked his interest in this classical result. The author also thanks Enrique Trevi\~no for useful discussions and feedback that improved the overall quality of this document.


\end{document}